\documentclass[10pt]{amsart}
\usepackage{latexsym,amscd,amssymb} 

\theoremstyle{plain}
  \newtheorem{theorem}{Theorem}[]
  \newtheorem{proposition}[theorem]{Proposition}
  \newtheorem{lemma}[theorem]{Lemma}

\theoremstyle{definition}

 \theoremstyle{remark}

\numberwithin{equation}{section}

\def\RR{{\mathbb R}}

\def\CC{{\mathbb C}}

\def\Hom{{\mathrm{Hom}}}
\def\Hilb{{\mathrm{Hilb}}}

\def\Sym{{\mathrm{Sym}}}
\def\rootht{{\mathrm{ht}}}

\begin{document}

\title[Fake degrees for reflection actions on roots]
{Fake degrees for reflection actions on roots}

\author{Victor Reiner}

\address{School of Mathematics\\
University of Minnesota\\
Minneapolis, MN 55455}
\email{reiner@math.umn.edu}

\author{Zhiwei Yun}

\address{Dept. of Mathematics\\
MIT\\
Cambridge, MA 02139}
\email{zyun@math.mit.edu}

\thanks{First, second authors partially supported by the 
NSF grants DMS-0601010, DMS-0969470.}



\begin{abstract}
A finite irreducible real reflection group of rank $\ell$ and Coxeter
number $h$ has root system of cardinality $h \cdot \ell$.
It is shown that the {\it fake degree} for the permutation action on 
its roots is divisible by $[h]_q=1+q+q^2+\cdots+q^{h-1}$, and that
in simply-laced types it equals
$
[h]_q \cdot \sum_{i=1}^\ell q^{d^*_i}
$
where $d^*_i=e_i-1$ are the {\it codegrees} and $e_i$ are the {\it exponents}.
\end{abstract}

\maketitle

\section{Introduction}
\label{fake-degrees-section}

Consider a complex reflection group $W \subset GL(V)$
with $V=\CC^\ell$, acting by linear substitutions on the polynomial algebra
$S=\Sym(V^*) \cong \CC[x_1,\ldots,x_n]$.  Both
Shephard and Todd \cite{ShephardTodd} and 
Chevalley \cite{Chevalley} proved 
that the invariant subalgebra is again a polynomial algebra  
$
S^W = \CC[f_1,\ldots,f_\ell]
$ 
for some homogeneous polynomials $f_i$, and
that the {\it coinvariant algebra} $S/I$ where
$I=(f_1,\ldots,f_\ell)$ carries a graded
version of the regular representation.
Thus for any finite-dimensional $\CC W$-module $U$, the intertwiner space
$
\Hom_W(U, S/I) \cong 
\left( U^* \otimes S/I \right)^W
$ 
is a graded $\CC$-vector space, whose {\it $q$-dimension}
or {\it Hilbert series} has been called its {\it fake degree}
$$
f^U(q) = \Hilb( \,\, \Hom_W( U, S/I ) \,\, ,  \,\, q ).
$$
Since 
$
f^U(1)= \dim_\CC \Hom(U, \CC W) =\dim_\CC U,
$ 
one may regard $f^U(q)$ as a $q$-analogue of the degree $\dim_\CC U$.
For example, the fake degree $f^{V^*}(q)$ of
the {\it dual reflection representation} $V^*$ is
determined by the {\it degrees} $d_1 \leq \cdots \leq d_{\ell}$ of the 
invariants $f_1,\ldots,f_\ell$ via\footnote{
This follows as a consequence of Solomon's result \cite{Solomon}
that the $W$-invariant {\it differential forms} with polynomial
coefficients $S \otimes \wedge^k V$ form a free $S^W$-module
with basis elements $df_{i_1} \wedge \cdots \wedge df_{i_k}$.
}
$
f^{V^*}(q) = \sum_{i=1}^\ell q^{d_i-1}.
$
One also defines the 
{\it codegrees} $d^*_1 \leq \cdots \leq  d^*_\ell$
via the fake degree polynomial 
$
f^{V}(q) = \sum_{i=1}^\ell q^{d^*_i+1}
$
of the representation $V$ itself.

We focus here on the case where $W$ acts 
on $V = \CC^\ell$ as the complexification of 
an irreducible {\it real} reflection group, so that one
has $V \cong V^*$ and $f^V(q)=f^{V^*}(q)$.
In this setting, one defines the 
{\it exponents} $(e_1,\ldots,e_\ell)$ by $e_i=d_i-1=d^*_i+1$,
and the {\it Coxeter number} $h=d_\ell$.
Choose a {\it root system} $\Phi$, containing
one opposite pair $\{\pm \alpha\}$ of normals
to each reflecting hyperplane, stable under the $W$-action.  
Given any $W$-stable subset $\Phi'$ of $\Phi$, we will
consider the fake degree polynomial $f^{\Phi'}(q):=f^U(q)$
for the $W$-permutation action $U=\CC \Phi'$.
Recall \cite[Chap. VI, \S 1, no. 2]{Bourbaki}, 
\cite[S 3.18]{Humphreys} that the 
cardinality of $\Phi$ has formula 
$
|\Phi|= h \cdot \ell.
$

\begin{theorem}
\label{main-result}
Let $W$ be an irreducible finite real reflection group, with root system $\Phi$,
and Coxeter number $h$.  Then for any $W$-stable subset of $\Phi'$ of $\Phi$,
\begin{enumerate}
\item[(i)]
$f^{\Phi'}(q)$ is divisible by $[h]_q=\frac{1-q^n}{1-q}$, and
\item[(ii)]
when $W$ is simply-laced, 
$
f^\Phi(q)=[h]_q \cdot ( q^{d^*_1} + \cdots +q^{d^*_\ell}).
$
\end{enumerate}
\end{theorem}

After posting this article to the arXiv, the authors learned that 
assertion (ii) of Theorem~\ref{main-result} appears 
in work of Stembridge \cite[Lemma 4.3(c,d)]{Stembridge2}, 
where it is proven via essentially the same method as in
Section~\ref{second-assertion-section} below.
Furthermore, Stembridge gives an explicit factorization 
of $f^{\Phi'}(q)/[h]_q$ in the general crystallographic case, where 
$\Phi'$ can be either of the two $W$-orbits of roots, short or long,
using a notion of {\it short exponents} for $W$.

\section{Proof of assertion $(i)$}

In the proof, one may assume without loss of generality
that $W$ acts transitively on the subset $\Phi'$ of $\Phi$.
The desired divisibility will then be deduced from
Lemma~\ref{coset-CSP-lemma} below,
applied to a {\it Coxeter element} of $W$.
The statement of the lemma involves Springer's notion \cite{Springer} of
a {\it regular element} $c$ in $W$, with a {\it regular eigenvalue} $\zeta$,
meaning $c(v) = \zeta v$ for an eigenvector $v$ 
lying on none of the reflecting hyperplanes for $W$.
Then $c$ and $\zeta$ have the same multiplicative order $n$.
Denote by $C$ the cyclic subgroup $\langle c \rangle$ generated by $c$.

\begin{lemma} \cite[Thm. 8.2]{RStantonWhite}
\label{coset-CSP-lemma}
Let $W$ be a complex reflection group acting transitively
on a finite set $X$, and $c$ in $W$ a regular element of 
order $n$, with a regular eigenvalue $\zeta$.
Then for all $m$, the fake degree $f^{X}(q):=f^U(q)$ for
the $W$-permutation action $U=\CC X$ satisfies
$$
 f^X(\zeta^m)  \quad =  \quad \#\{ x \in X: c^m(x)=x\}.
$$
In particular, $f^X(q)$ is divisible by $[n]_q$
if and only if $C$ acts freely on $X$.

\end{lemma}
\begin{proof}
For the sake of completeness, we recall the
proof from \cite{RStantonWhite}.
Springer \cite{Springer}
extended the work of Shephard-Todd and Chevalley
by proving one has an isomorphism $W \times C$-representations 
\begin{equation}
\label{Springer-isomorphism}
S/I \cong \CC W
\end{equation}
where $W$ acts as before, and where $C$ acts on
$\CC W$ via {\it right-translation}, and on
$S/I$  via scalar substitutions $c(x_i) = \zeta^{-1}\cdot x_i$.
Equivalently, $c$ scales the $d^{th}$ homogeneous component $(S/I)_d$
by the scalar $\zeta^{-d}$.

Now identify the transitive $W$-permutation representation $\CC X$
with a coset action $\CC[W/W']$ for some subgroup $W'$ of $W$.
Then one has an isomorphism 
$\Hom_W( \CC [W/W'], S/I) \cong (S/I)^{W'}$, and one can
reformulate the fake degree:
\begin{equation}
\label{fake-degree-as-invariant-hilb}
f^{X}(q) = \Hilb( (S/I)^{W'} , q ).
\end{equation}
Taking $W'$-fixed spaces in \eqref{Springer-isomorphism}
give an isomorphism of $C$-representations 
\begin{equation}
\label{fixed-space-isomorphism}
(S/I)^{W'} \cong (\CC W)^{W'} \cong \CC X
\end{equation}
and the result now follows by comparing the
trace of $c^m$ on the two ends of \eqref{fixed-space-isomorphism}.
\end{proof}

\noindent
To finish the proof of assertion (i),
one applies Lemma~\ref{coset-CSP-lemma} to a finite {\it real} reflection 
group $W$, with {\it Coxeter generators} $S=\{s_1,\ldots,s_\ell\}$, and 
with $c=s_1 s_2 \cdots s_\ell$ a {\it Coxeter element}.
It is known that all Coxeter elements lie in a single $W$-conjugacy class,
that they have multiplicative order $h=d_\ell$, and that they are regular
elements having $\zeta=e^{\frac{2 \pi i}{h}}$ as a regular eigenvalue; see
\cite[\S 3.16, 3.17]{Humphreys}.  
Furthermore, it is known \cite[Chap. V, \S 1, no. 11]{Bourbaki} 
that the cyclic group $C$ generated by a Coxeter element $c$ 
acts freely on the roots $\Phi$.  Assertion (i) now follows from 
Lemma~\ref{coset-CSP-lemma}.

\section{Proof of assertion $(ii)$}
\label{second-assertion-section}

We first recall a bit more of the
root geometry for finite real reflection groups, in order to
further reformulate the fake degree $f^{\Phi'}(q)$;
see e.g. \cite[Chapters 1, 5]{Humphreys}.

Assume $W$ is the complexification of a real reflection
group acting on $V_\RR \cong \RR^\ell$, that preserves a positive
definite inner product $(-,-)$ on $V_\RR$. 
The reflecting hyperplanes
dissect $V_\RR$ into open simplicial 
cones called {\it chambers}, which are permuted simply-transitively
by $W$.  Choosing one such chamber $C$ to be the
{\it dominant chamber}, every $W$-orbit contains exactly
one point in its closure $\bar{C}$.  The root system
decomposes as $\Phi=\Phi_+ \sqcup -\Phi_+$,
where the positive roots $\Phi_+$ are those
having positive inner product with the points of $C$.
This also distinguishes the subset 
of {\it simple roots} $\{\alpha_1,\ldots,\alpha_\ell\}$
inside $\Phi_+$, whose nonnegative linear combinations contain $\Phi_+$,
and whose corresponding {\it simple reflections} 
$S=\{s_1,\ldots,s_\ell\}$ gives rise to a Coxeter presentation
$(W,S)$ for $W$.  
The above discussion implies that every 
$W$-orbit of roots contains a unique {\it dominant}
representative $\alpha_0$ lying in $\bar{C}$, whose isotropy
subgroup $W_{\alpha_0}$ 
is a {\it standard parabolic subgroup} generated by some\footnote{Although
we will not need this information here, the table
at the beginning of Section~\ref{remarks-section} lists the
type for these standard parabolic subgroups $W_{\alpha_0}$.
When $W$ is crystallographic and $\alpha_0$ is the highest root, 
$W_{\alpha_0}$ is generated by the simple reflections of $W$ not adjacent 
to the extra node $s_{0}$ in the extended Dynkin diagram for the 
affine Weyl group $\widetilde{W}$.} subset $S$.

\begin{proposition}
\label{highest-root-reformulation-prop}
Let $W$ be a finite real reflection group $W$ with root system $\Phi$
and positive roots $\Phi_+$.
Let $\Phi'$ be a $W$-orbit of roots, with unique dominant
representative $\alpha_0$.  Then the fake degree for the $W$-permutation
action on $\Phi'$ can be expressed as
$$
f^{\Phi'}(q)
= \sum_{\alpha \in \Phi'} q^{d(\alpha_0,\alpha)}
$$
where $d(\alpha_0,\alpha)$ is the Coxeter group length $\ell_S(w)$
of the minimum length representative $w$ for the
coset $wW_{\alpha_0}=\{u \in W: u(\alpha_0)=\alpha\}$.
\end{proposition}
\begin{proof}
Note that $S$ is a free $S^W$-module, because $S^W=\CC[f_1,\ldots,f_\ell]$
is a polynomial ring.  One obtains $S^W$-module splittings for the ring inclusions 
$S^{W_{\alpha_{0}}} \subset S$ and $S^W \subset S^{W_{\alpha_{0}}}$ by averaging over 
$W_{\alpha_0}$ and over coset representatives for $W/W_{\alpha_0}$, respectively.
Hence $S^{W_{\alpha_0}}$ is also a free $S^W$-module, with
$$
f^{\Phi'}(q) = \Hilb( (S/I)^{W_{\alpha_0}}, q) 
=\frac{ \Hilb( S^{W_{\alpha_0}},q) }{\Hilb(S^W,q)}.
$$
For any standard parabolic subgroup $W'$ of $W$, such as $W'=W_{\alpha_0}$ or
$W'=W$ itself, one has \cite[\S 3.15]{Humphreys} that
$
\Hilb(S^{W'},q)^{-1}= (1-q)^{\ell} \sum_{w \in W'} q^{\ell_S(w)}.
$
Therefore
\begin{equation}
\label{min-coset-rep-formula}
f^{\Phi'}(q) = \frac{ \sum_{w \in W} q^{\ell_S(w)} }
                   { \sum_{w \in W_{\alpha_0}} q^{\ell_S(w)} }
 = \sum_{w} q^{\ell_S(w)}
\end{equation}
where in this last sum, $w$ runs over the minimum-length coset representatives
for the cosets $wW_{\alpha_0}$ in $W/W_{\alpha_0}$.
\end{proof}

The crux of the proof of assertion (ii) will be the
following lemma\footnote{This lemma is similar in spirit
to results of Stembridge \cite[\S2,3]{Stembridge1} on a quantity that he
calls the {\it depth} $d(\alpha)$ of the root $\alpha$, closely related to the
quantity $d(\alpha_0,\alpha)$ defined here.}.  It 
relates, for simply-laced root systems with
highest root $\alpha_0$, the quantity $d(\alpha_0,\alpha)$ to
the {\it root height} of $\alpha$, which we recall here; 
see \cite[Chap. VI, \S8]{Bourbaki}, \cite[\S 3.20]{Humphreys},
\cite[\S3]{Stembridge1} for further discussion.
When $W$ is a crystallographic root system $\Phi$, with 
simple roots $\{\alpha_1,\ldots,\alpha_\ell\}$, for
every root $\alpha$ in $\Phi$ the
unique expression $\alpha=\sum_{i=1}^{\ell} c_i \alpha_i$
has {\it integer} coefficients $c_i$, and one
defines the {\it height} $\rootht(\alpha)=\sum_{i=1}^\ell c_i$.
There is a unique {\it highest root} $\alpha_0$,
achieving the maximum height $\rootht(\alpha_0)=h-1$, and this
highest root $\alpha_0$ is always dominant.

\begin{lemma}
\label{M-V-lemma}
Let $W$ be a simply-laced root Weyl group with root system
$\Phi$, positive roots $\Phi_+$, 
and highest root $\alpha_0$.  Then any root $\alpha$ in $\Phi$ has
$$
d(\alpha_0,\alpha)=
\begin{cases}
\rootht(\alpha_0)-\rootht(\alpha) & \text{ if } \alpha \in \Phi_+,\\
\rootht(\alpha_0)-\rootht(\alpha)-1 & \text{ if } \alpha \in -\Phi_+.\\
\end{cases}
$$
\end{lemma}
\begin{proof}
Rescale all roots $\alpha$ so that $(\alpha,\alpha)=2$,
and consequently $(\alpha,\beta)$ lies in $\{0,\pm 1,\pm 2\}$ for all pairs
of roots $\alpha, \beta$.  For any simple root, the 
formula 
$$
s_i(\beta)=\beta-(\beta,\alpha_i) \alpha_i
$$
shows that applying the simple reflection $s_i$ to a root 
$\beta \neq \pm \alpha_i$ has the following effect on its height:
$$
\rootht(s_i \beta) =
\begin{cases}
\rootht(\beta) & \text{ if } (\beta,\alpha_i)=0 \\
\rootht(\beta)+1 & \text{ if } (\beta,\alpha_i)=-1 \\
\rootht(\beta)-1 & \text{ if } (\beta,\alpha_i)=+1. 
\end{cases}
$$
When $\beta = \pm \alpha_i$, one has $\rootht(\beta)=\pm 1$,
and $\rootht(s_i(\beta))=-\rootht(\beta)=\mp 1$.

Consequently, when starting with the highest root $\alpha_0$, and applying
a sequence of simple reflections $s_i$, the height can drop by at most one
at each stage, except when one crosses from a simple root to its negative.
This implies that the expression on the right side in the lemma
(call it $b(\alpha)$) gives a lower bound on the length $\ell_S(w)$ for any $w$ sending
$\alpha_0$ to $\alpha$.  Thus $d(\alpha_0,\alpha) \geq b(\alpha)$.  

To show $d(\alpha_0,\alpha) \leq b(\alpha)$, induct on $b(\alpha)$.
In the base case $b(\alpha)=0$, so $\alpha=\alpha_0$ and $d(\alpha_0,\alpha)=0$ also.
In the inductive step, $b(\alpha) \neq 0$ implies $\alpha \neq \alpha_0$,
so (as we are in the simply-laced case)
$\alpha$ is not dominant, and there exists some simple root
$\alpha_i$ with $(\alpha,\alpha_i)<0$.
It suffices to show that $b(s_i \alpha) = b(\alpha)-1$.

If $(\alpha,\alpha_i)=-1$ then $\rootht(s_i\alpha)=\rootht(\alpha)+1$,
and either both $\alpha, s_i(\alpha)$ lie in $\Phi_+$ or both lie in $-\Phi_+$, 
so $b(s_i \alpha) = b(\alpha)-1$.  

If $(\alpha,\alpha_i)=-2$ then $\alpha=-\alpha_i$, so that 
$s_i\alpha=+\alpha_i$, and again $b(s_i\alpha)=b(\alpha)-1$.
\end{proof}

The proof of assertion (ii) requires one more
well-known fact \cite[\S 3.20]{Humphreys},
relating the distribution of root heights to the exponents $e_i = d^*_i+1$:
\begin{equation}
\label{exponent-height-relation}
\sum_{\alpha \in \Phi_{+}} q^{\rootht(\alpha)} 
= \sum_{i=1}^\ell (q^1+q^2+\cdots+q^{e_i}).
\end{equation}
For $W$ simply-laced, there is only one orbit $\Phi$,
whose dominant root $\alpha_0$ is the highest root, with $\rootht(\alpha_0)=h-1$.
Combining Proposition~\ref{highest-root-reformulation-prop},
Lemma~\ref{M-V-lemma}, \eqref{exponent-height-relation} gives
$$
\begin{aligned}
f^{\Phi}(q)
&= \sum_{\alpha \in \Phi_+} q^{h-1-\rootht(\alpha)}
+\sum_{\alpha \in -\Phi_+} q^{h-2-\rootht(\alpha)}\\
&=\sum_{i=1}^{\ell} 
(q^{h-e_i-1}+q^{h-e_i}+\cdots +q^{h-2})+
(q^{h-1}+q^h+\cdots+q^{h+e_i-2}) \\
&= (1-q)^{-1} \sum_{i=1}^{\ell} (q^{h-e_i-1}-q^{h+e_i-1}) \\
&= (1-q)^{-1} \sum_{i=1}^{\ell} (q^{e_i-1}-q^{h+e_i-1})
\end{aligned}
$$
where the last equality used the fact \cite[\S 3.16]{Humphreys} that $h-e_i= e_{\ell+1-i}$.
Therefore
$$
f^{\Phi}(q)
= \frac{1-q^h}{1-q} \cdot \sum_{i=1}^{\ell} q^{e_i-1}\\
= [h]_q \cdot \sum_{i=1}^{\ell} q^{d^*_i}
$$
as desired.

\section{Remarks and questions}
\label{remarks-section}

\subsection{Further divisibilities}

The table below tabulates the polynomial $f^{\Phi'}(q)/[h]_q$ for 
root orbits $\Phi'$ in all real reflection groups.  In the crystallographic
types $A-E$, this can also be deduced from Stembridge's {\it exponent data} 
\cite[Table 4.1]{Stembridge2} together with 
his factorization \cite[Lemma 4.2(c,d)]{Stembridge2}.  The last column tabulates
the additional data $\gcd( [h]_q, \sum_{i=1}^\ell q^{d^*_i})$, relevant for 
Proposition~\ref{no-common-factors} below.  
\vskip.1in
\begin{tabular}{|c||c|c|c|c|c|}
\hline
$W$ & $h$ & $\Phi'=W.\alpha_0$  & $W_{\alpha_0}$ type & $f^{\Phi'}(q)/[h]_q$ & $\gcd([h]_q, \sum_i q^{d^*_i})$ \\ \hline\hline\hline
$A_{n-1}$ & $n$ & $\Phi$ & $A_{n-3}$ & $[n-1]_q$ & $1$\\ \hline \hline
$B_n$ & $2n$ & $\{\pm e_i \pm e_j\}$ & $A_1 \times B_{n-2}$ & ${[n-1]_{q^2}}$ &
$[n]_{q^2}$ \\   
      &      & $\{\pm e_i\}$& $B_{n-1}$ & $1$ & \\ \hline\hline
$D_n$ & $2(n-1)$ & $\Phi$ & $A_1 \times D_{n-2}$ 
   & $\frac{[n-2]_{q^2} [n]_q}{ [2]_q}$ & $1$ \\ \hline \hline
$E_6$ & $12$ & $\Phi$ & $A_5$ & $[2]_{q^4} [3]_{q^3}$ & $1$\\ \hline
$E_7$ & $18$ & $\Phi$ & $D_6$ & $\frac{[2]_{q^6}}{[2]_{q^2}} [7]_{q^2}$ & $1$ \\ \hline
$E_8$ & $30$ & $\Phi$ & $E_7$ & $[2]_{q^{10}} [4]_{q^6}$ & $1$ \\ \hline \hline
$F_4$ & $12$ & either orbit & $B_3$ & $[2]_{q^4}$ & $[2]_{q^6}$ \\ \hline\hline
$H_3$ & $10$ & $\Phi$& $A_1 \times A_1$ & $[3]_{q^2}$ & $1$ \\ \hline
$H_4$ & $30$ & $\Phi$& $H_3$ & $[2]_{q^6}[2]_{q^{10}}$ & $1$ \\ \hline\hline
$I_2(m)$ & $m$ & either orbit & $A_1$ & $1$ &  
    $1$ if $\frac{m}{2}$ odd \\ 
$m$ even & & & & & $[2]_{q^2}$ if $\frac{m}{2}$ even \\ \hline
$I_2(m)$ & $m$ & $\Phi$ & $-$ & $[2]_q$ & $1$ \\
$m$ odd & & & & & \\ \hline
\end{tabular}

\vskip.1in
\noindent
The table exhibits case-by-case two facts for which we lack uniform proofs.

\begin{proposition}
\label{no-common-factors}
For finite real $W$ with one root orbit,
$\gcd([h]_q,\sum_{i=1}^{\ell} q^{d^*_i})=1$. 
\end{proposition}
\noindent
Using \eqref{Springer-isomorphism},
Proposition~\ref{no-common-factors} is equivalent to the assertion that,
when $W$ has only one orbit of roots, every power $c^m$ of a 
Coxeter element $c$ acts on $V$ with nonzero trace.

\begin{proposition}
\label{fake-deg-divisibility}
For finite real $W$ which are {\bf at most doubly-laced},
meaning that its Coxeter presentation relations $(s_is_j)^{m_{ij}}=e$ all have 
$m_{ij} \leq 4$, every $W$-stable root subset $\Phi'$ has 
fake degree $f^{\Phi'}(q)$ divisible by $\sum_{i=1}^{\ell} q^{d^*_i}$.
\end{proposition}

\subsection{Original motivation}
We originally observed Theorem~\ref{main-result} case-by-case
while computing the fake degree of a certain {\it irreducible} 
representation of simply-laced $W$, arising
naturally in \cite[Chapter 3]{RSaliolaWelker}.  
One can decompose the $W$-permutation representation 
$\CC[\Phi']$ of any real reflection group $W$ 
on a root orbit $\Phi'$ into two
direct summands, namely its {\it symmetric} and 
{\it antisymmetric} components $\CC[\Phi']^{+}, \CC[\Phi']^{-}$
with respect to the $W$-equivariant involution that 
simultaneously swaps each $+\alpha, -\alpha$.
A straightforward calculation then shows the following.

\begin{proposition}
\label{three-related-reps-prop}
Let $W$ be a finite real reflection group $W$ and $\Phi'$ an orbit
of its roots.  Then any one of the three 
fake degrees for $\CC[\Phi'], \CC[\Phi']^{+}, \CC[\Phi']^{-}$
determines the others via the relations
$
f^{\Phi'}(q) = f^{\Phi',^+}(q) + f^{\Phi',-}(q)
$
and 
$
f^{\Phi',-}(q) =q \cdot f^{\Phi',+}(q).
$
\end{proposition}

\noindent
It was further shown in 
\cite[Chapter 3]{RSaliolaWelker} that,
for irreducible real reflection groups $W$, and any root orbit $\Phi'$,
the antisymmetric component $\CC[\Phi']^-$
has $W$-irreducible decomposition which is {\it multiplicity-free}.
In the simply-laced case, it has only two irreducible
constituents:  
$
\CC[\Phi]^-= V \oplus U
$
where $V$ is the reflection representation $V$ of degree $\ell$,
and $U$ is another $W$-irreducible, of degree $|\Phi^+|-\ell=\frac{h-2}{2} \cdot \ell$.
Using Proposition~\ref{three-related-reps-prop},
one can check that Theorem~\ref{main-result}(ii) is equivalent to
the assertion that this $W$-irreducible $U$ has fake degree
$
f^U(q)= q^2 \cdot \frac{[h-2]_q}{[2]_q} \cdot \sum_{i=1}^{\ell} q^{e_i}.
$

\subsection{M-V cycles}
Lemma~\ref{M-V-lemma} has a geometric interpretation.
It is well-known that for a standard parabolic subgroup $W'$ of a Weyl group
$W$ associated to simple complex algebraic group $G$ and Borel subgroup $B$, 
one can identify the invariant subalgebra $(S/I)^{W'}$ 
with the cohomology $H^*(G/P)$ of $G/P$ where $P=\langle B, W'\rangle$.  
The Schubert cell decomposition of $G/P$ lets one
express its Poincar\'e polynomial in terms of
lengths of minimal coset representatives for $W/W'$.
The expression \eqref{min-coset-rep-formula} then arises
in this way when $W'=W_{\alpha_0}$ for a dominant root $\alpha_0$.

When $\alpha_0$ happens to be the highest root of a simply-laced root
system, the cone over the variety $G/P$
also arises as a Schubert variety in the affine Grassmannian. The cell decomposition of $G/P$ as above can be used to give a decomposition of this cone into {\it Mirkovi\'c-Vilonen cycles} introduced in \cite{MirkovicVilonen}. In this
picture, the dimension formula for the Mirkovi\'c-Vilonen cycles is
equivalent to Lemma~\ref{M-V-lemma};
see  Mirkovi\'c and Vilonen \cite[Theorem 3.2]{MirkovicVilonen} with $\lambda = \alpha_0$,
and also Ng\^{o} and Polo \cite[Lemme 7.4]{NgoPolo}.

\subsection{A-D-E quivers?}
For simply-laced $W$, the $W$-action permuting the roots
can be modeled by {\it reflection functors} acting on the
the bounded derived category of quiver representations,
with a Coxeter element $c$ corresponding to the
{\it Auslander-Reiten translation}.  Here the $W$-equivariant
map from an object to its dimension vector 
factors through the quotient category that
mods out by the square of the shift map; see the
discussion of the {\it periodic Auslander-Reiten quiver} by 
Kirillov and Thind \cite{KirillovThind}.
Does Theorem~\ref{main-result}(ii) reflect something
lurking in this quiver picture?

\end{document}